\documentstyle[amscd,amssymb,verbatim,epsf]{amsart}

\begin{document}
\theoremstyle{plain}
\newtheorem{Thm}{Theorem}
\newtheorem{Cor}{Corollary}
\newtheorem{Ex}{Example}
\newtheorem{Con}{Conjecture}
\newtheorem{Main}{Main Theorem}
\newtheorem{Lem}{Lemma}
\newtheorem{Prop}{Proposition}

\theoremstyle{definition}
\newtheorem{Def}{Definition}
\newtheorem{Note}{Note}

\theoremstyle{remark}
\newtheorem{notation}{Notation}
\renewcommand{\thenotation}{}

\errorcontextlines=0
\numberwithin{equation}{section}
\renewcommand{\rm}{\normalshape}%

\title[Lagrangian Curves]%
   {Lagrangian curves on spectral curves of monopoles}

\author{Brendan Guilfoyle}
\address{Brendan Guilfoyle\\
          Department of Mathematics\\
          Institute of Technology, Tralee\\
          Clash\\
          Tralee\\
          Co. Kerry\\
          Ireland.}
\email{brendan.guilfoyle@@ittralee.ie}
\author{Madeeha Khalid}
\address{Madeeha Khalid\\
          Department of Mathematics\\
          Institute of Technology, Tralee\\
          Clash\\
          Tralee\\
          Co. Kerry\\
          Ireland.}
\email{madeeha.khalid@@staff.ittralee.ie}

\author{Jos\'e J. Ram\'on Mar\'i}
\address{Jos\'e J. Ram\'on Mar\'i\\
          IT Tralee Research Institute\\
          Institute of Technology, Tralee\\
          Clash\\
          Tralee\\
          Co. Kerry\\
          Ireland.}
\email{jose.ramon.mari@@research.ittalee.ie}

\keywords{Neutral Kaehler structure, monopoles, Lagrangian curves}
\subjclass{Primary: 53A25; Secondary: 81T13}
\date{18th October, 2007}

\begin{abstract}
We study Lagrangian points on smooth holomorphic curves in T${\mathbb P}^1$ equipped with a natural neutral
K\"ahler structure, and prove that they must
form real curves. By virtue of the identification of T${\mathbb P}^1$ with the space ${\mathbb L}({\mathbb E}^3)$ 
of oriented affine lines in Euclidean 3-space ${\mathbb E}^3$, these Lagrangian curves give rise to ruled surfaces in
${\mathbb E}^3$, which we prove have zero Gauss curvature. 

Each ruled surface is shown to be the tangent lines
to a curve in ${\mathbb E}^3$, called the edge of regression of the ruled surface.  We give an alternative characterization
of these curves as the points  in ${\mathbb E}^3$ where the number of oriented lines in the complex curve $\Sigma$ 
that pass through the point is less than the degree of $\Sigma$. We then apply these results to the spectral curves of 
certain monopoles and construct the ruled surfaces and edges of regression generated by the Lagrangian curves.
\end{abstract}

\maketitle

The space ${\mathbb L}({\mathbb E}^3)$ of oriented affine lines in Euclidean 3-space ${\mathbb E}^3$ can be identified
with the total space TS$^2$ of the tangent bundle to the 2-sphere. The standard complex structure on ${\mathbb P}^1$
lifts to a complex structure ${\mathbb J}$ on T${\mathbb P}^1$, and hence ${\mathbb L}({\mathbb E}^3)$.
This is well-known and has found a variety of uses, most notably
in the construction of monopoles in ${\mathbb E}^3$ \cite{hitch1}. What is less well-known is the canonical symplectic 
structure on ${\mathbb L}({\mathbb E}^3)$ which is compatible with ${\mathbb J}$ and enjoys many remarkable
geometric properties \cite{gak4} \cite{gak5} \cite{salvai}.

Together the complex and symplectic structures form a K\"ahler structure with the property that the associated
metric is of signature ($++--$). In this paper we consider Lagrangian points on complex curves in T${\mathbb P}^1$ 
equipped with this neutral K\"ahler structure. At such points the metric induced on the surface vanishes.

The only complex curves that are Lagrangian at every point
are the oriented normals to planes and spheres in ${\mathbb E}^3$, and we exclude these curves. 
Our first main result is:

\vspace{0.1in}

\noindent{\bf Main Theorem 1}:

{\it
Let $\Sigma$ be a smooth compact complex curve in T ${\mathbb P}^1$.

\begin{enumerate}
\item[(i)] The branch points of the composition $\Sigma\hookrightarrow {\rm T}\;{\mathbb P}^1\rightarrow{\mathbb{P}}^1$ 
           are Lagrangian,
\item[(ii)] there do not exist any isolated Lagrangian points on $\Sigma$,
\item[(iii)] if $C\subset\Sigma$ is a Lagrangian curve, then the associated ruled surface in ${\mathbb{E}}^3$ has 
              zero Gauss curvature.
\end{enumerate}
}

\vspace{0.1in}

We show that a ruled surfaces in part (iii) is tangent to a curve in ${\mathbb E}^3$, called the edge of regression
of the ruled surface. We can characterize this curve in ${\mathbb E}^3$ another way. Given a holomorphic curve
$\Sigma$ consider the number of oriented lines in $\Sigma$ that pass through a given point in ${\mathbb E}^3$. 
Generically, we show that this number is the degree of $\Sigma$ when it is considered as a curve in ${\mathbb{P}}^3$.
Moreover, we prove:

\vspace{0.1in}

\noindent{\bf Main Theorem 2}:

{\it
Let $\Sigma$ be a smooth complex curve in T${\mathbb P}^1$. Then the genus of $\Sigma$ is $(m-1)^2$ for $m=1,2,3...$ 
and a generic point in ${\mathbb E}^3$ has $2m$ distinct oriented lines of $\Sigma$ passing through it. 
This is the maximum number of distinct oriented lines of $\Sigma$ that can pass through a point 
(the minimum number being one).

The points in ${\mathbb E}^3$ lying on less than $2m$  distinct oriented lines of $\Sigma$  form the edges of regression 
of the ruled surfaces generated by the Lagrangian curves on $\Sigma$.
}

\vspace{0.1in}

We apply this to the spectral curves of the charge 2 and charge 3 monopoles 
\cite{HMR}, paying particular attention to the ruled surfaces and the edges of regression.

We show that the Lagrangian curves in the charge 2 case consists of 4 disjoint curves, 2 passing through the branch
points and 2 which do not. These curves generate two ruled surfaces in ${\mathbb E}^3$ which have edges of regression 
an ellipse in the $x^1x^2-$plane, and a hyperbola in the $x^1x^3-$plane. The eccentricity of the ellipse is
 $k$, while that of the hyperbola is with $1/k$, where $k$ is the spectral parameter of the monopole. 

For the charge 3 monopole we show that the set of Lagrangian points consists of 10 (intersecting) simple closed curves. 
These 
fall naturally into two classes: 6 of the curves are rulings of planes, in particular, 6 planes that form
the edges of a tetrahedron in ${\mathbb E}^3$. The remaining 4 curves are rulings of flat, non-planar surfaces forming the 
faces of the tetrahedron.

Finally, in the
last section, we discuss the results and place them in the more general setting of null curves on neutral K\"ahler
surfaces.  

\section{Lagrangian Curves on Holomorphic Curves}

The total space TS$^2$ of the tangent bundle to the 2-sphere is a 4-manifold with a natural complex structure
defined as follows. Let $\xi$ be the standard holomorphic coordinate on ${\mathbb P}^1$ given by stereographic 
projection from the south pole.
Let $\eta$ be the corresponding complex coordinate in the fibre of the bundle $\pi: {\mbox {TS}}^2\rightarrow 
{\mbox S}^2$ 
obtained by identifying the pair of complex numbers ($\xi,\eta$) with the tangent vector
\[
\eta\frac{\partial}{\partial \xi}+\bar{\eta}\frac{\partial}{\partial \bar{\xi}}\in {\mbox T}_\xi {\mbox S}^2.
\]
This yields holomorphic coordinates on TS$^2-\pi^{-1}\{\mbox{south pole}\}$, which can be supplemented by an analogous
coordinate system on TS$^2-\pi^{-1}\{\mbox{south pole}\}$. The transition functions on the overlap are holomorphic
and so we obtain a complex structure ${\mathbb J}$ on TS$^2$. For short we write T${\mathbb{P}}^1$ for TS$^2$ with 
this complex structure.

The space ${\mathbb L}({\mathbb E}^3)$ of oriented affine lines in Euclidean 3-space ${\mathbb E}^3$ can be identified
with the total space TS$^2$, and thus inherits the above complex structure. In fact,
this complex structure is natural in the sense that it is invariant under the action of the Euclidean group acting on 
${\mathbb L}({\mathbb E}^3)$ \cite{gak4}.

Geometric data can be transferred between ${\mathbb{L}}({\mathbb{E}}^3)$ and ${\mathbb{E}}^3$ by use of a 
correspondence space.

\begin{Def}
The map $\Phi:{\mathbb{L}}({\mathbb{E}}^3)\times{\mathbb{R}}\rightarrow{\mathbb{E}}^3$ is defined to take 
$(\gamma,r)\in{\mathbb{L}}({\mathbb{E}}^3)\times{\mathbb{R}}$ to the point  in ${\mathbb{E}}^3$ on the oriented line 
$\gamma$ that lies a distance $r$ from the point on the line closest to the origin.
\end{Def}

\newpage

The double fibration below gives us the correspondence between the points in ${\mathbb{L}}({\mathbb{E}}^3)$ 
and oriented lines in ${\mathbb{E}}^3$:
we identify a point $\gamma$ in ${\mathbb{L}}({\mathbb{E}}^3)$ with $\Phi\circ\pi_1^{-1}(\gamma)\subset{\mathbb{E}}^3$, which is an oriented line.
Similarly, a point $p$ in ${\mathbb{E}}^3$ is identified with the 2-sphere $\pi_1\circ\Phi^{-1}(p)\subset{\mathbb{L}}({\mathbb{E}}^3)$, which consists of all of the 
oriented lines through the point $p$.

\vspace{0.7in}

\unitlength0.5cm

\begin{picture}(16,6)
\put(8.1,5.8){$\pi_1$}
\put(8.0,7.5){${\mathbb{L}}({\mathbb{E}}^3)\times{\mathbb{R}}$}
\put(10,7){\vector(1,-1){2.8}}
\put(11.5,5.8){$\Phi$}
\put(8.4,3){${\mathbb{L}}({\mathbb{E}}^3)$}
\put(9,7){\vector(0,-1){3}}
\put(13,3.4){${\mathbb{E}}^3$}
\end{picture}

If $\Phi((\xi,\eta),r)=(z(\xi,\eta,r),t(\xi,\eta,r))$, then it has the following coordinate expression \cite{gak2}:
\begin{equation}\label{e:coord}
z=\frac{2(\eta-\bar{\eta}\xi^2)+2\xi(1+\xi\bar{\xi})r}{(1+\xi\bar{\xi})^2}
\qquad\qquad
t=\frac{-2(\eta\bar{\xi}+\bar{\eta}\xi)+(1-\xi^2\bar{\xi}^2)r}{(1+\xi\bar{\xi})^2},
\end{equation}
where $z=x^1+ix^2$, $t=x^3$ and ($x^1$, $x^2$, $x^3$) are Euclidean
coordinates in ${\mathbb{E}}^3$.

The symplectic structure on T${\mathbb{P}}^1$ can be motivated in a number of ways - for our purposes we will
simply use its coordinate expression (more details can be found in \cite{gak4}):
\begin{equation}
\Omega=\frac{2}{(1+\xi\bar{\xi})^2}\left(
  d\eta \wedge d\bar{\xi}+d\bar{\eta} \wedge d\xi
   +\frac{2(\xi\bar{\eta}-\bar{\xi}\eta)}{1+\xi\bar{\xi}}d\xi\wedge d\bar{\xi}
\right).
\end{equation}
This is clearly a closed non-degenerate 2-form on T${\mathbb P}^1$ which is compatible with ${\mathbb J}$
\[
\Omega({\mathbb J}({\mathbb X}),{\mathbb J}({\mathbb Y}))=\Omega({\mathbb X},{\mathbb Y})
\qquad\qquad
{\mbox {for all }}\quad {\mathbb X},{\mathbb Y}\in \mbox{T}_\gamma \mbox{T}{\mathbb P}^1.
\] 
The symplectic 2-form is also invariant under the action of the Euclidean group \cite{gak4}.

The neutral K\"ahler metric on T${\mathbb P}^1$ is defined by ${\mathbb{G}}(\cdot,\cdot)=\Omega({\mathbb J}\cdot,\cdot)$
and has local coordinate expression
\[
{\mathbb{G}}=\frac{2i}{(1+\xi\bar{\xi})^2}\left(
  d\eta d\bar{\xi}-d\bar{\eta} d\xi
   +\frac{2(\xi\bar{\eta}-\bar{\xi}\eta)}{1+\xi\bar{\xi}}d\xi d\bar{\xi}
\right).
\]

\vspace{0.2in}

We now consider a real 2-dimensional surface $\Sigma\subset\mbox{T}{\mathbb P}^1$. 

\begin{Def}
A point $\gamma\in\Sigma$ is a {\it complex point} if 
${\mathbb J}:$T$_\gamma{\mathbb P}^1\rightarrow$T$_\gamma{\mathbb P}^1$.
A surface $\Sigma$ is called a {\it complex curve} (or {\it holomorphic curve}) if every point of $\Sigma$ is a 
complex point.

A point $\gamma\in\Sigma$ is said to be a {\it Lagrangian point} if the symplectic 2-form $\Omega$ pulled back
to T$_\gamma\Sigma$ is zero. A surface $\Sigma$ is called {\it Lagrangian} if every point of 
$\Sigma$ is a Lagrangian point.

\end{Def}

The only real surfaces in T${\mathbb P}^1$ that are both complex and Lagrangian at every point are the oriented normal 
lines to planes and spheres in ${\mathbb E}^3$. In what follows we exclude this case.

We now prove:

\vspace{0.2in}

\begin{Main} \label{t:mt1}

Let $\Sigma$ be a smooth compact complex curve in T ${\mathbb P}^1$.

\begin{enumerate}
\item[(i)] The branch points of the composition $\Sigma\hookrightarrow {\rm T}\;{\mathbb P}^1\rightarrow{\mathbb{P}}^1$ 
           are Lagrangian,
\item[(ii)] there do not exist any isolated Lagrangian points on $\Sigma$,
\item[(iii)] if $C\subset\Sigma$ is a Lagrangian curve, then the associated ruled surface in ${\mathbb{E}}^3$ has 
              zero Gauss curvature.
\end{enumerate}
\end{Main}

\begin{pf}
\vspace{0.1in}

Let $\Sigma$ be a (real) surface in T${\mathbb P}^1$. About any point $\gamma\in\Sigma$ there is a local parameterization 
${\mathbb C}\rightarrow {\mathbb C}^2$:$\nu\mapsto(\xi(\nu,\bar{\nu}),\eta(\nu,\bar{\nu}))$, where we assume, 
without loss of generality, that $\gamma$ does not lie in the fibre over the south pole.

The real surface $\Sigma$ is holomorphic iff about each point of $\Sigma$ we have $\sigma=0$ where \cite{gak2}
\begin{equation}\label{e:comp}
\sigma=\partial \xi\bar{\partial}\eta-\bar{\partial} \xi\partial\eta,
\end{equation}
$\partial$ being differentiation with respect to the parameter $\nu$. 

On the other hand, by pulling back the 2-form $\Omega$, we see that a real surface is Lagrangian at a point 
$\gamma\in\Sigma$ iff we have $\lambda={\mathbb I}{\mbox m}\;\rho\;$=0 at $\gamma$, where 
\begin{equation}\label{e:rho}
\rho=\partial\eta\bar{\partial}\bar{\xi}-\bar{\partial}\eta\partial\bar{\xi}-\frac{2\bar{\xi}\eta}{1+\xi\bar{\xi}}
\left(\partial\xi\bar{\partial}\bar{\xi}-\bar{\partial}\xi\partial\bar{\xi}\right).
\end{equation}

It is also clear that $\Sigma$ is locally the graph of a section of the canonical bundle 
$\pi:{\mbox T}{\mathbb P}^1\rightarrow{\mathbb P}^1$ iff
\[
\partial \xi\bar{\partial}\bar{\xi}-\bar{\partial} \xi\partial\bar{\xi}\neq0.
\]

We now turn to the proofs of statements (i) to (iii).  Let $\Sigma$ be a holomorphic curve in T${\mathbb P}^1$ so that 
$\sigma=0$. 

\vspace{0.1in}
\noindent{\bf Proof of (i)}:

Suppose that $\gamma\in\Sigma$ is a branch point. Then, as it is smooth, the curve osculates the fibre of the 
canonical bundle at $\gamma$ and so
\[
\partial \xi\bar{\partial}\bar{\xi}-\bar{\partial} \xi\partial\bar{\xi}=0,
\]
at $\gamma$.

A short calculation shows that
\[
\rho\bar{\rho}-\sigma\bar{\sigma}=(\partial \eta\bar{\partial}\bar{\eta}-\bar{\partial} \eta\partial\bar{\eta})
            (\partial \xi\bar{\partial}\bar{\xi}-\bar{\partial} \xi\partial\bar{\xi}),
\]
which therefore vanishes at $\gamma$. However, $\sigma=0$ and so we conclude that at a branch point $\rho=0$. 
In particular, $\lambda={\mathbb I}{\mbox m}\;\rho=0$ at $\gamma$, and so the point is Lagrangian, as claimed.

\vspace{0.1in}

\noindent{\bf Proof of (ii)}:

We argue by contradiction. Let $\gamma\in\Sigma$ be an isolated Lagrangian point, which we assume, 
without loss of generality, does not lie on $\pi^{-1}\{{\mbox {south pole}}\}$.
Thus there exists an open neighbourhood U$\subset\Sigma$ containing $\gamma$ such that 
\[
\lambda(\gamma)=0 \qquad\qquad \lambda\left(U-\{\gamma\}\right)\neq 0.
\]  

First suppose that the Lagrangian point $\gamma$ is a branch point of the mapping 
$\Sigma\hookrightarrow {\mbox {T}}{\mathbb P}^1\rightarrow{\mathbb{P}}^1$. 
Then, since the projection restricted to $\Sigma$ is not of maximal rank, $\Sigma$ cannot be locally parameterized by a 
section of this bundle. However, as $\Sigma$ is smooth, we can use the fibre coordinate as a local parameter about 
$\gamma$:
$\eta\rightarrow(\xi(\eta,\bar{\eta}),\eta)$. Since $\sigma=0$, by equation (\ref{e:comp}) we have that 
$\bar{\partial}\xi=0$ and by equation (\ref{e:rho})
\[
\lambda={\textstyle{\frac{i}{2}}}(\partial\xi-\bar{\partial}\bar{\xi}).
\] 
Thus $\lambda$ is the imaginary part of a holomorphic function and therefore its zeros cannot be isolated. This 
proves the claim when $\gamma$ is a branch point.

Now suppose that the Lagrangian point $\gamma$ is not a branch point. Then we can parameterize a neighbourhood U of
$\gamma$ on $\Sigma$ by a local section of the bundle: $\xi\rightarrow(\xi,\eta(\xi,\bar{\xi}))$. Since $\Sigma$
is holomorphic, by equation (\ref{e:comp}) we have $\bar{\partial}\eta=0$ and by (\ref{e:rho}) 
\[
\lambda={\mathbb I}{\mbox{m}}\;(1+\xi\bar{\xi})^2\partial\left(\frac{\eta}{(1+\xi\bar{\xi})^2}\right).
\]
A short computation shows then that $\lambda$ satisfies the second order equation
\[
\partial\bar{\partial}\lambda+\frac{2\lambda}{1+\xi\bar{\xi}}=0.
\]
The strong maximum and minimum principles ({\it cf.} Theorem 2.2 of \cite{GaT}) imply that zeros of $\lambda$ 
cannot be isolated.
For, suppose that $\lambda\geq 0$ on U. Then $\lambda$ is superharmonic on U: $\partial\bar{\partial}\lambda\leq0$, 
while $\lambda(\gamma)=\inf_U\lambda$. By the strong minimum principle, $\lambda$ is constant, in fact zero, on U, which 
is a contradiction (as we have ruled out the oriented normals to planes and spheres).

If $\lambda\leq 0$ on U, then $\lambda$ is subharmonic on U: $\partial\bar{\partial}\lambda\geq0$, 
while $\lambda(\gamma)=\sup_U\lambda$. By the strong maximum principle, $\lambda$ is constant, in fact zero, on U, which 
again is a contradiction.

We conclude that none of the Lagrangian points on the holomorphic curve can be isolated.

\vspace{0.1in}

\noindent{\bf Proof of (iii)}:
 
Classically, a {\it ruled surface} is a 1-parameter family of oriented lines in ${\Bbb{E}}^3$.
From our point of view, a ruled surface is a real curve $C$ in ${\mbox {T}}{\mathbb P}^1$. Suppose that
this curve is given locally by $s\rightarrow (\xi(s),\eta(s))$. Then, by equations (\ref{e:coord}), the ruled surface is
\[
z(r,s)=\frac{2[\eta(s)-\bar{\eta}(s)\xi(s)^2]+2\xi(s)[1+\xi(s)\bar{\xi}(s)]r}{[1+\xi(s)\bar{\xi}(s)]^2},
\]
\[
t(r,s)=\frac{-2[\eta(s)\bar{\xi}(s)+\bar{\eta}(s)\xi(s)]+[1-\xi(s)^2\bar{\xi}(s)^2]r}{[1+\xi(s)\bar{\xi}(s)]^2},
\]
where, as before, $z=x^1+ix^2$, $t=x^3$ and ($x^1$, $x^2$, $x^3$) are Euclidean
coordinates in ${\mathbb{E}}^3$, and $r$ is an affine parameter along the lines of the ruling. 

By a straightforward, if lengthy, curvature calculation, the Gauss curvature of such a ruled surface is found to be:
\begin{equation}\label{e:gauss}
K=-\frac{(1+\xi\bar{\xi})^2\left[{\mathbb{I}}{\mbox {m}}\left((1+\xi\bar{\xi})\dot{\eta}\dot{\bar{\xi}}
                   +2\xi\bar{\eta}\dot{\xi}\dot{\bar{\xi}}\right)\right]^2}
{\left|(1+\xi\bar{\xi})\dot{\eta}-2\bar{\xi}\eta\dot{\xi}+(1+\xi\bar{\xi})r\dot{\xi}\right|^4},
\end{equation}
where a dot represents differentiation with respect to $s$.

If the real curve $C$ lies on a holomorphic curve we have $\dot{\eta}=\partial\eta\dot{\xi}$ and the Gauss
curvature simplifies to
\[
K=-\frac{\lambda^2}
    {(\lambda+(r+\psi)^2)^2},
\]
where $\rho=\psi+i\lambda$ as in (\ref{e:rho}). Along a Lagrangian curve $\lambda=0$, and so the Gauss curvature of the 
ruled surface vanishes for such a curve.

This completes the proof of the theorem.

\end{pf}

\vspace{0.2in}

Flat ruled surfaces (referred to as {\it developable} ruled surfaces) were studied in classical surface theory. Aside from 
rulings of
a plane, other examples of flat ruled surfaces include generalized cones and cylinders. In fact:

\begin{Thm}
The ruled surface generated by a Lagrangian curve on a holomorphic curve is the tangent lines to an oriented curve in 
${\mathbb{E}}^3$.
\end{Thm}
\begin{pf}

A well-known result of classical surface theory states that every developable surface can be subdivided into 
portions of a cylinder, a cone or the tangent line to a curve in ${\mathbb{E}}^3$ (for example, see Thm. 58.3 of 
\cite{krey}).
For the situation stated in the theorem we eliminate the first two possibilities: generalised cylinder and cone, as 
follows. 

The generalised cylinder is obtained by translating an oriented line along a curve 
in the plane orthogonal to the oriented line. Clearly
the direction of the lines in this ruling do not change and such a curve in T${\mathbb P}^1$ must
therefore lie in a fibre
of the projection $\pi:{\mbox {T}}{\mathbb P}^1\rightarrow {\mathbb P}^1$. This is not the case for Lagrangian 
curves on holomorphic curves (as we have ruled out the case
of the oriented normals to planes) and so the ruled surface cannot be a  generalised cylinder. 

The generalised cone is a 1-parameter family of oriented lines passing through a fixed point p in ${\mathbb E}^3$. 
Such a curve in T${\mathbb P}^1$ lies on the holomorphic sphere of all oriented lines through p. Thus, were a Lagrangian 
curve $C$ on a holomorphic curve to form a generalized cone, it would lie on the intersection of two holomorphic curves,
an impossibility for a 1-dimensional set.
 
\end{pf}

Thus every Lagrangian curve on a holomorphic curve gives rise to a curve in ${\mathbb E}^3$: the 
{\it edge of regression} of the ruled surface. This subset of ${\mathbb E}^3$ can be defined in a different way, as we
show in the next section.

\section{A Different Characterization}

In this section we give an alternative characterization of the set of Lagrangian points on a smooth complex curve
in T${\mathbb P}^1$. In particular, we prove:

\vspace{0.2in}
\begin{Main}
Let $\Sigma$ be a smooth complex curve in T${\mathbb P}^1$. Then the genus of $\Sigma$ is $(m-1)^2$ for $m=1,2,3...$ 
and a generic point in ${\mathbb E}^3$ has $2m$ distinct oriented lines of $\Sigma$ passing through it. This is 
the maximum number of distinct oriented lines of $\Sigma$ that can pass through a point (the minimum number being one).

The points in ${\mathbb E}^3$ lying on less than $2m$  distinct oriented lines of $\Sigma$  form the edges of regression 
of the ruled surfaces generated by the Lagrangian curves on $\Sigma$.
\end{Main}
\begin{pf}
The set of oriented lines through a point in ${\mathbb{E}}^3$ forms a global holomorphic section of the 
complex vector bundle
T${\mathbb P}^1\rightarrow{\mathbb P}^1$. In particular, an oriented line $(\xi,\eta)$ passes through a fixed point 
$(z,t)\in{\mathbb{C}}\oplus{\mathbb{R}}={\mathbb{E}}^3$ iff
\begin{equation}\label{e:laghol}
\eta={\textstyle{\frac{1}{2}}}\left(z-2t\xi-\bar{z}\xi^2\right).
\end{equation}

Let $F_2={\mathbb P}({\mathcal O}\oplus {\mathcal O}(2))$ be the Hirzebruch surface obtained by 
the quotient of ${\mathbb C}^2- \{0\}\times{\mathbb C}^2-\{0\}$ by the equivalence relation 
\[
(x_0,x_1,y_0,y_1)\sim(ax_0,ax_1,by_0,ba^2y_1),
\]
for $a,b\in{\mathbb C}^*$.
Projecting onto the first factor, we see that $F_2$ is, in fact, a ${\mathbb P}^1$ bundle over ${\mathbb P}^1$. 

T${\mathbb P}^1$ is isomorphic to $F_2-E_\infty$, where $E_\infty=\{(x_0,x_1,0,y_1)\}/\sim$ is the infinity section. 
That is, $F_2$ is isomorphic to 
T${\mathbb P}^1$ after the one point compactification of each of the fibres of the canonical projection
T${\mathbb P}^1\rightarrow{\mathbb P}^1$.

Now, the Picard group of $F_2$ is $\mbox{Pic}(F_2)={\mathbb Z}[h]\oplus{\mathbb Z}[f]$ where $[h]$ and $[f]$ are 
the divisor classes of the holomorphic
sections and fibres, respectively, of the bundle $F_2\rightarrow{\mathbb P}^1$. These classes have intersection
pairing
\[
[h]\cdot[h]=2, \qquad\qquad [f]\cdot[f]=0, \qquad\qquad [h]\cdot[f]=1, 
\] 
and in this basis
\begin{equation}\label{e:einfty}
[E_\infty]=[h]-2[f].
\end{equation}
To see this, note that $[E_\infty]=k[h]+l[f]$ for some $k,l\in{\mathbb Z}$, and taking the intersection with 
$[h]$ and $[f]$ we find that, since $[E_\infty]\cdot[h]=0$ and $[E_\infty]\cdot[f]=1$, $2k+l=0$ and $k=1$, 
which yield (\ref{e:einfty}).

Thus, if $\Sigma$ is a compact complex curve in T${\mathbb P}^1=F_2-E_\infty$, we have
\[
0=[E_\infty]\cdot[\Sigma]=[h]\cdot[\Sigma]-2[f]\cdot[\Sigma].
\]
We conclude that $[h]\cdot[\Sigma]=2[f]\cdot[\Sigma]=2m$ for some $m\in{\mathbb N}$, and therefore $\Sigma$ 
intersects a generic holomorphic
section of T${\mathbb P}^1$ in $2m$ distinct points. This extends to generic holomorphic sections of the form 
(\ref{e:laghol})
and so we conclude that a generic point in ${\mathbb E}^3$ has $2m$ distinct oriented lines of $\Sigma$ passing through it.

To determine the genus $g$ of $\Sigma$, we use the adjunction formula
\[
\left.\left(K_{F_2}+[\Sigma]\right)\right|_\Sigma=K_\Sigma,
\]
where $K_{F_2}$ and $K_\Sigma$ are the canonical bundles of $F_2$ and $\Sigma$, respectively. 
This implies that
\begin{equation}\label{e:adj}
[\Sigma]\cdot[\Sigma]+K_{F_2}\cdot[\Sigma]=2g-2.
\end{equation}
By Lemma V 2.10 of \cite{hart} and equation (\ref{e:einfty})
\[
K_{F_2}=-2[E_\infty]-4[f]=-2[h],
\]
and a routine calculation shows that $[\Sigma]=m[h]$ for some $m\in{\mathbb N}$. Thus
\[
[\Sigma]\cdot[\Sigma]+K_{F_2}\cdot[\Sigma]=2m^2-4m,
\]
and by (\ref{e:adj}) we find that $g=(m-1)^2$ as claimed.

Away from the branch points, $\Sigma$ is given locally
by a holomorphic section $\eta=F(\xi)$. Let p$\in{\mathbb E}^3$ and consider the oriented lines in $\Sigma$ that pass 
through p. If p has coordinates $(z,t)$ as above, then
we are seeking to find the roots of
\[
G=F-{\textstyle{\frac{1}{2}}}\left(z-2t\xi-\bar{z}\xi^2\right).
\]

An oriented line $(\xi_0,F(\xi_0))$ is a multiple root of $G$ iff $G(\xi_0)=0$ and $\partial G(\xi_0)=0$.
These are equivalent to
\begin{equation}\label{e:multroot1}
F(\xi_0)={\textstyle{\frac{1}{2}}}\left(z-2t\xi_0-\bar{z}\xi_0^2\right),
\end{equation}
\begin{equation}\label{e:multroot2}
(1+\xi\bar{\xi})^2\partial\left.\left(\frac{F}{(1+\xi\bar{\xi})^2}\right)\right|_{\xi_0}
       =-\frac{z\bar{\xi}_0+\bar{z}\xi_0+t(1-\xi_0\bar{\xi}_0)}{1+\xi_0\bar{\xi}_0}.
\end{equation}

Suppose now that $(\xi_0,F(\xi_0))\in\Sigma$ is a multiple root of $G$ and so the above equations hold for 
some ($z,t$). Then clearly
\[
{\mathbb I}{\mbox{m}}\left.\left[(1+\xi\bar{\xi})^2\partial\left(\frac{F}{(1+\xi\bar{\xi})^2}\right)
\right]\right|_{\xi_0}=0,
\]
and so the point is Lagrangian.

Conversely, suppose that $(\xi_0,F(\xi_0))\in\Sigma$ is Lagrangian, and define
\begin{equation}\label{e:caustic}
r_0=-(1+\xi\bar{\xi})^2\partial\left.\left(\frac{F}{(1+\xi\bar{\xi})^2}\right)\right|_{\xi_0}.
\end{equation}
By the Lagrangian condition this is a real number. Now consider the following point $(z,t)\in{\mathbb E}^3$:
\[
z=\frac{2[F(\xi_0)-\overline{F(\xi_0)}\xi_0^2]+2\xi_0(1+\xi_0\bar{\xi}_0)r_0}{(1+\xi_0\bar{\xi}_0)^2},
\]
\[
t=\frac{-2[F(\xi_0)\bar{\xi}_0+\overline{F(\xi_0)}\xi_0]+(1-\xi_0^2\bar{\xi}_0^2)r_0}{(1+\xi_0\bar{\xi}_0)^2},
\]
A calculation shows that equations (\ref{e:multroot1}) and (\ref{e:multroot2}) hold, so that less than $d$ 
distinct oriented lines pass through the point 
($z,t$), and this point lies on the oriented line $(\xi_0,F(\xi_0))\in\Sigma$.

\end{pf}

\begin{Note}
$F_2$ is a resolution of the quadric cone $Q$ in ${\mathbb P}^3$, where $E_\infty$ is the exceptional divisor. 
The intersection number $2m$ above is the degree of $\Sigma$ when considered as a curve in ${\mathbb P}^3$.
This can be described in local coordinates as follows. 

Let [$z_0:z_1:z_2:z_3$] be homogenous coordinates 
on ${\mathbb{P}}^3$. As described above, T${\mathbb P}^1$ can be 
identified with an open subset of $F_2$, and $F_2$ can be mapped to ${\mathbb{P}}^3$ by 
\[
(x_0,x_1,y_0,y_1)\mapsto[x_0^2y_0:x_0x_1y_0:x_1^2y_0:y_1]
\]
Clearly, $F_2$ maps to the 
quadric cone $Q$ given by $z_0z_2=z_1^2$, and the infinity section $E_\infty$ maps to the vertex
$p=[0:0:0:1]$ of the cone. 

The result is an identification of T${\mathbb P}^1$ with $Q-p$ which can be written locally 
\begin{equation}\label{e:quadcone}
(\xi,\eta)\longleftrightarrow[1:\xi:\xi^2:\eta].
\end{equation}

A curve $\Sigma$ in T${\mathbb P}^1$ maps to a curve on $Q-p\subset {\mathbb P}^3$. But for any 
complex curve in ${\mathbb P}^3$ there is a well-defined degree 
\[
d={\mbox{deg }}\Sigma=\#(\Sigma\cap{\mathbb{P}}^2)
\]
where ${\mathbb{P}}^{2}$ is a generic complex plane in ${\mathbb{P}}^3$.

Now, a generic ${\mathbb{P}}^2$ is given by
\[
a_1z_1+a_2z_2+a_3z_3+a_4z_4=0,
\]
so that, using the identification (\ref{e:quadcone}),  $Q\cap{\mathbb{P}}^2$ is
\[
a_1+a_2\xi+a_3\xi^2+a_4\eta=0.
\]
For $a_4\neq0$ this defines a global holomorphic section of T${\mathbb P}^1\rightarrow{\mathbb P}^1$ and so $d$ is 
exactly the total number of points of intersection of $\Sigma$ and a generic holomorphic section, that is, $d=2m$.
\end{Note}

\begin{Note}
Main Theorem 2 restricts the genus of a smooth holomorphic curve in T${\mathbb P}^1$. In contrast, taking the
oriented normal lines to a smooth surface of genus $g$ in ${\mathbb E}^3$ we obtain a smooth Lagrangian surface
of genus $g$ in T${\mathbb P}^1$, and so there exist smooth Lagrangian surfaces of any genus in T${\mathbb P}^1$.
\end{Note}

\vspace{0.2in}

\section{Lagrangian Curves on Spectral Curves}

The complex structure on T${\mathbb P}^1$ plays a crucial role in Hitchin's approach to BPS monopoles
for SU(2) Yang-Mills-Higgs theory in ${\mathbb E}^3$ \cite{hitch1}. Each such monopole can be constructed from its spectral curve: a 
compact holomorphic curve in T${\mathbb P}^1$ which is given in our coordinates by
\[
\eta^m+\alpha_1(\xi)\eta^{m-1}+...+\alpha_m(\xi)=0,
\]
where each $\alpha_j$ is a complex polynomial of degree less than or equal to $2j$ satisfying the reality condition
\[
\alpha_j(\xi)=(-1)^j\xi^{2j}\overline{\alpha_j\left(-\frac{1}{\bar{\xi}}\right)}.
\]
Here $m\in{\mathbb N}$ is the charge of the monopole, and the moduli space of gauge inequivalent charge $m$ monopoles is a 
(4$m$-1)-dimensional manifold. 

We now consider the Lagrangian curves on the spectral curves of the
charge 2 and the tetrahedrally symmetric charge 3 monopole.

\subsection{The Charge 2 Monopole}

The moduli space of gauge inequivalent charge 2 monopoles is a 
7-dimensional manifold. However, 6 of these degrees of freedom can be removed by the action of the Euclidean
group, so that there is only a 1-parameter family of charge 2 monopoles.
The spectral curve $\Sigma$ can be described in local coordinates on T${\mathbb{P}}^1$ by 
\cite{ManSut}
\begin{equation}\label{e:ch2}
\eta^2=\alpha\left[k^2(1+\xi^4)-2(2-k^2)\xi^2 \right],
\end{equation}
where the spectral parameter $k$ satisfies $0\leq k <1$ and 
\[
\alpha=\left[{\textstyle \frac{1}{2}}\int_0^{\scriptstyle \frac{\pi}{2}}\frac{d\theta}{\sqrt{1-k^2\sin^2\theta}} \right]^2.
\]
The parameter $k$ measures the separation of the two monopoles. As $k\rightarrow 1$, 
the separation of the monopoles goes to infinity. The case $k=0$, which corresponds to the two monopoles 
coinciding at the origin, leads to a singular spectral curve and will be excluded in what follows. 

For $k\neq0$ the spectral curve is a smooth compact complex curve of genus 1 in T${\mathbb P}^1$ which
double covers ${\mathbb P}^1$ with 4 branch points 
$\gamma_1,\gamma_2,\gamma_3,\gamma_4$. These points lie at the 4 real points
\[
\xi=\pm\sqrt{\frac{2-k^2\pm2\sqrt{1-k^2}}{k^2}}.
\]
We start by finding where the symplectic form on T${\mathbb{P}}^1$ pulled back to $\Sigma$ vanishes:
\[
\Omega|_\Sigma=2\left[\partial\left(\frac{\eta}{(1+\xi\bar{\xi})^2}\right)
    -\bar{\partial}\left(\frac{\bar{\eta}}{(1+\xi\bar{\xi})^2}\right)\right]d\xi\wedge d\bar{\xi}
    =\frac{4\lambda i}{(1+\xi\bar{\xi})^2}d\xi\wedge d\bar{\xi}=0.
\]
Differentiating (\ref{e:ch2}) we compute that
\[
\partial\eta=\frac{2\alpha}{\eta}\left[k^2\xi^3-(2-k^2)\xi\right],
\]
and so
\[
\partial\eta-\frac{2\bar{\xi}\eta}{1+\xi\bar{\xi}}=\frac{2\alpha}{\eta(1+\xi\bar{\xi})}\left[-k^2\bar{\xi}-(2-k^2)\xi
     +(2-k^2)\xi^2\bar{\xi}+k^2\xi^3\right].
\]

Upon squaring, the Lagrangian condition $\lambda=0$ becomes
\[
\frac{\left[-k^2\bar{\xi}-(2-k^2)\xi+(2-k^2)\xi^2\bar{\xi}+k^2\xi^3\right]^2}{k^2(1+\xi^4)-2(2-k^2)\xi^2}=
\frac{\left[-k^2\xi-(2-k^2)\bar{\xi}+(2-k^2)\xi\bar{\xi}^2+k^2\bar{\xi}^3\right]^2}{k^2(1+\bar{\xi}^4)-2(2-k^2)\bar{\xi}^2}.
\]
This can be simplified and factorized to
\[
k^2(1-k^2)(1-\xi\bar{\xi})(1+\xi\bar{\xi})^3(\xi-\bar{\xi})(\xi+\bar{\xi})=0.
\]

\vspace{0.1in}
\setlength{\epsfxsize}{4.5in}
\begin{center}
   \mbox{\epsfbox{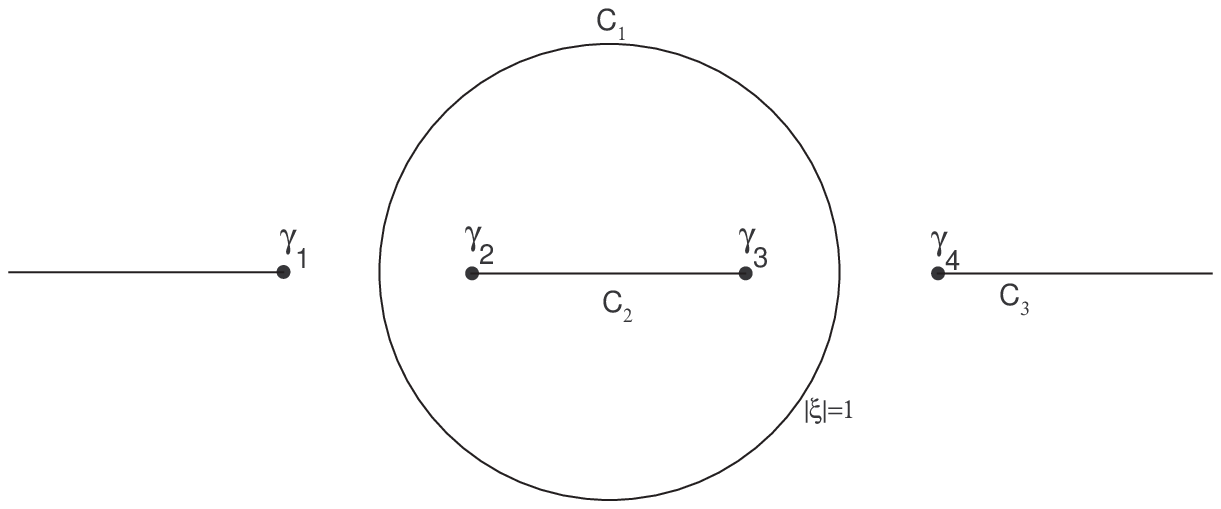}}
\end{center}
\vspace{0.1in}

Aside from the limits $k=0$ and $k=1$, we therefore have $|\xi|=1$, $\xi=\bar{\xi}$ or $\xi=-\bar{\xi}$. The last of these
does not lift to a Lagrangian curve on $\Sigma$, being an artefact of the squaring of the Lagrangian condition. Thus the 
Lagrangian points on $\Sigma$ project down to two curves on ${\mathbb P}^1$: the equator and a line of longitude, which
we now consider in detail.

Firstly, parameterize the equator C$_1$: $|\xi|=1$ by $\xi=e^{i\theta}$. Then, we compute
\[
\eta^2=\alpha\left[k^2(1+e^{4i\theta})-2(2-k^2)e^{2i\theta}\right]=-4\alpha(1-k^2\cos^2\theta)e^{2i\theta}.
\]
The curve C$_1$ lifts to two curves $\eta=\pm 2i\sqrt{\alpha}\sqrt{1-k^2\cos^2\theta}\;e^{i\theta}$ on $\Sigma$.
Substituting this in equation (\ref{e:coord}) we find the ruled surface in ${\mathbb E}^3$ is
\[
x^1=\mp2\sqrt{\alpha}\sqrt{1-k^2\cos^2\theta}\;\sin\theta+r\cos\theta
\qquad
x^2=\pm2\sqrt{\alpha}\sqrt{1-k^2\cos^2\theta}\;\cos\theta+r\sin\theta,
\]
\[
x^3=0.
\]

The other Lagrangian points are given by $\xi=\bar{\xi}$. Suppose $\xi=\tan(\theta/2)$ for $-\pi\leq\theta\leq\pi$
then
\[
\eta^2=\frac{\alpha(k^2-\sin^2\theta)}{\cos^4(\theta/2)} 
\qquad
{\mbox {or }}
\qquad
\eta=\pm\frac{\sqrt{\alpha}\sqrt{k^2-\sin^2\theta}}{\cos^2(\theta/2)}. 
\]
To be Lagrangian, clearly we must have $\eta\in{\mathbb R}$, so $\sin^2\theta\leq k^2$. This implies that
$\theta$ must lie in the following domains:
\[
-\pi\leq\theta\leq-\pi+\sin^{-1}k
\qquad
-\sin^{-1}k\leq\theta\leq\sin^{-1}k
\qquad
\pi-\sin^{-1}k\leq\theta\leq\pi.
\]
The first and last of these domains are connected: they form the curve C$_3$ when projected onto ${\mathbb P}^1$ and 
contain the branch points $\gamma_1$ and $\gamma_4$. The
middle domain forms the curve C$_2$ which passes through the branch points $\gamma_2$ and $\gamma_3$. 
Of course each interval lifts to two copies in T${\mathbb P}^1$ joined at the branch points, forming circles in the 
total space.

The ruled surfaces in ${\mathbb E}^3$, obtained by inserting our parameterized curves in equation (\ref{e:coord}), are
\[
x^1=\mp2\sqrt{\alpha}\sqrt{k^2-\sin^2\theta}\cos\theta+r\sin\theta 
\qquad x^2=0,
\]
\[
x^3=\pm2\sqrt{\alpha}\sqrt{k^2-\sin^2\theta}\sin\theta+r\cos\theta .
\]

We thus have shown that the Lagrangian curves on the charge 2 spectral curve yields rulings of the $x^1x^2-$ and
$x^1x^3-$planes. We can find the edge of regression of such a ruling by using the procedure given by (\ref{e:caustic}). 
In particular, the edge of regression of these lines is obtained by substituting
\[
r=-\partial\eta+\frac{2\bar{\xi}\eta}{1+\xi\bar{\xi}},
\]
in equation (\ref{e:coord}). The result for the lift of curve C$_1$ is 
\[
x^1=\mp\frac{2\sqrt{\alpha}\sin\theta}{\sqrt{1-k^2\cos^2\theta}}
\qquad
x^2=\pm\frac{2\sqrt{\alpha}(1-k^2)\cos\theta}{\sqrt{1-k^2\cos^2\theta}}
\qquad
x^3=0,
\]
This is an ellipse of eccentricity $k$, as can be seen by noting that the parameterized curve satisfies:
\[
\frac{(x^1)^2}{4\alpha}+\frac{(x^2)^2}{4\alpha(1-k^2)}=1.
\]
For the lift of curves C$_2$ and C$_3$ we obtain the edge of regression
\[
x^1=\mp\frac{2\sqrt{\alpha}k^2\cos\theta}{\sqrt{k^2-\sin^2\theta}}
\qquad x^2=0
\qquad
x^3=\pm\frac{2\sqrt{\alpha}(k^2-1)\sin\theta}{\sqrt{k^2-\sin^2\theta}},
\]
which is a hyperbola with eccentricity $1/k$. To see this, note that
\[
\frac{(x^1)^2}{4\alpha k^2}-\frac{(x^2)^2}{4\alpha(1-k^2)}=1.
\]

In the diagram below we show the ruled surface generated by the Lagrangian curves on the charge 2 monopole with $k=0.8$ 
- the edges of regression (a hyperbola and an ellipse) are clearly identifiable. Note that the 4 branch points of the
spectral curve are the two asymptotes of the hyperbola (counted once with each orientation).

\vspace{0.1in}
\setlength{\epsfxsize}{4.5in}
\begin{center}
   \mbox{\epsfbox{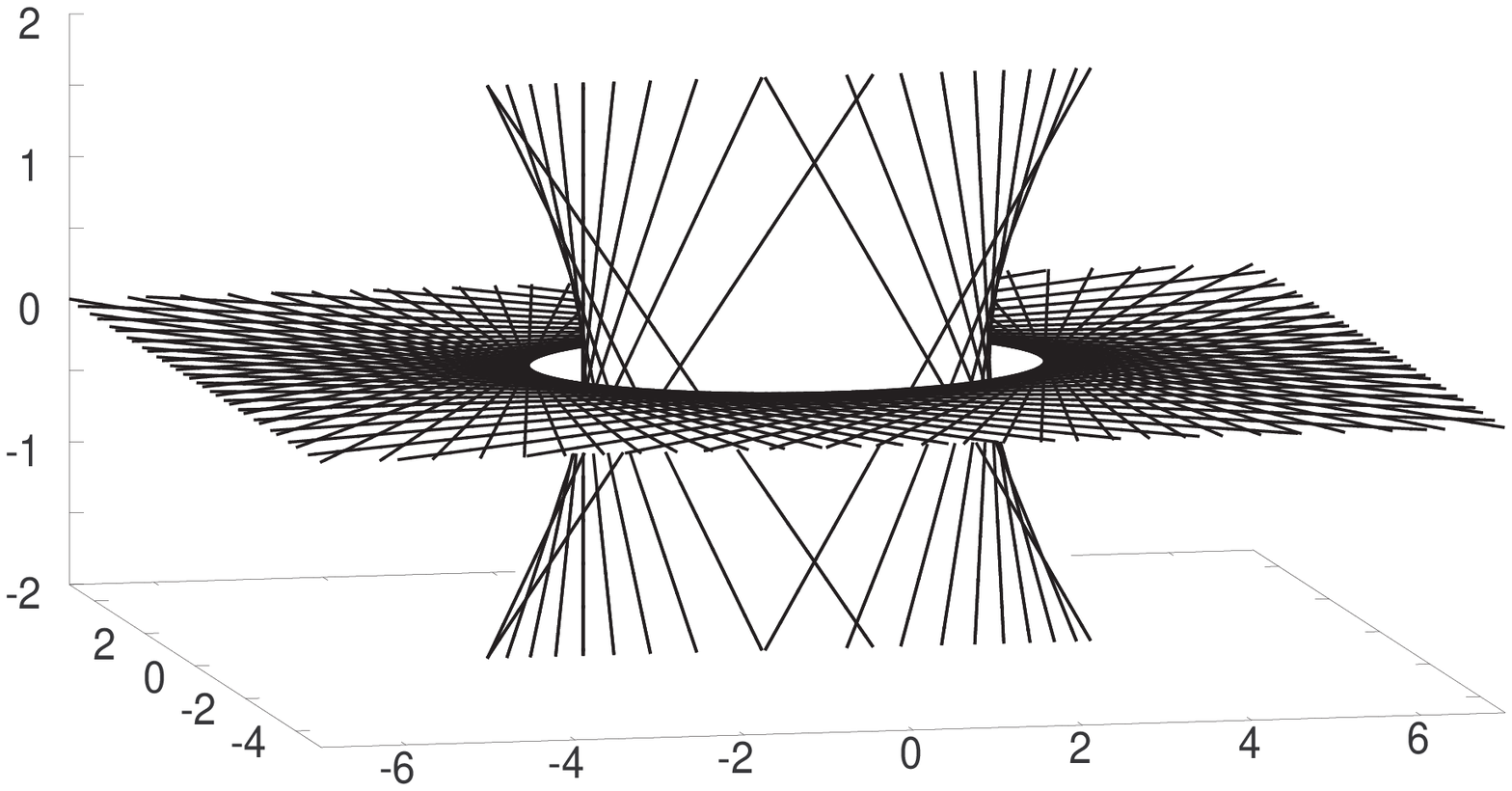}}
\end{center}
\vspace{0.1in}

\subsection{The Charge 3 Monopole}

The spectral curve $\Sigma$ of the tetrahedrally symmetric charge 3 monopole is \cite{HMR}:
\begin{equation}\label{e:ch3}
\eta^3=\frac{\Gamma{\textstyle{\left(\frac{1}{3}\right)}}^9}{48\sqrt{6}\pi^3}(1-5\sqrt{2}\xi^3-\xi^6).
\end{equation}

We start by pulling back the symplectic 2-form to $\Sigma$ and finding that the Lagrangian curves are given by
the equation:
\[
{\mathbb I}{\mbox m}\left[\eta^2\left(2\xi+5\sqrt{2}\bar{\xi}^2-5\sqrt{2}\xi\bar{\xi}^3+2\bar{\xi}^5\right)\right]=0.
\]
Substitute (\ref{e:ch3}) in this and cube the resulting equation to get
\[
{\mathbb I}{\mbox m}\left[\left(1-5\sqrt{2}\xi^3-\xi^6\right)^2\left(2\xi+5\sqrt{2}\bar{\xi}^2-5\sqrt{2}\xi\bar{\xi}^3+2\bar{\xi}^5\right)^3\right]=0.
\]
We factorize this to
\begin{equation}\label{e:lag_curve}
(\xi-\bar{\xi})(1+\xi\bar{\xi})^3(\xi^2+\xi\bar{\xi}+\bar{\xi}^2){\mathbb R}{\mbox {e}}(f(\xi,\bar{\xi}))=0,
\end{equation}
where
\begin{align}
f(\xi,\bar{\xi})=&-4 \xi^{9} \bar{\xi}^{9}-165 \sqrt{2} \xi^{9} \bar{\xi}^{6}-42 \xi^{9} \bar{\xi}^{3}+10 \sqrt{2} \xi^{9}
+162 \xi^{8} \bar{\xi}^{8}+810 \sqrt{2} \xi^{8} \bar{\xi}^{5}\nonumber\\
&-162 \xi^{8} \bar{\xi}^{2}-162 \xi^{7} \bar{\xi}^{7}-810 \sqrt{2} \xi^{7} \bar{\xi}^{4}
+162 \xi^{7} \bar{\xi}-2988 \xi^{6} \bar{\xi}^{6}+360 \sqrt{2} \xi^{6} \bar{\xi}^{3}\nonumber\\
&+42 \xi^{6}+8100 \xi^{5} \bar{\xi}^{5}-810 \sqrt{2} \xi^{5} \bar{\xi}^{2}-8100 \xi^{4} \bar{\xi}^{4}+810 \sqrt{2} \xi^{4} \bar{\xi}\nonumber\\
&+2988 \xi^{3} \bar{\xi}^{3}-165 \sqrt{2} \xi^{3}+162 \xi^{2} \bar{\xi}^{2}-162 \xi \bar{\xi}
+4\nonumber.
\end{align}

The equation (\ref{e:lag_curve}) defines the set of points on ${\mathbb {P}}^1$ obtained by the projection of the 
Lagrangian curves on $\Sigma$. The set can be plotted numerically, which we do below (after stereographic projection onto
the plane) on two scales:

\vspace{0.1in}
\setlength{\epsfxsize}{5.0in}
\begin{center}
   \mbox{\epsfbox{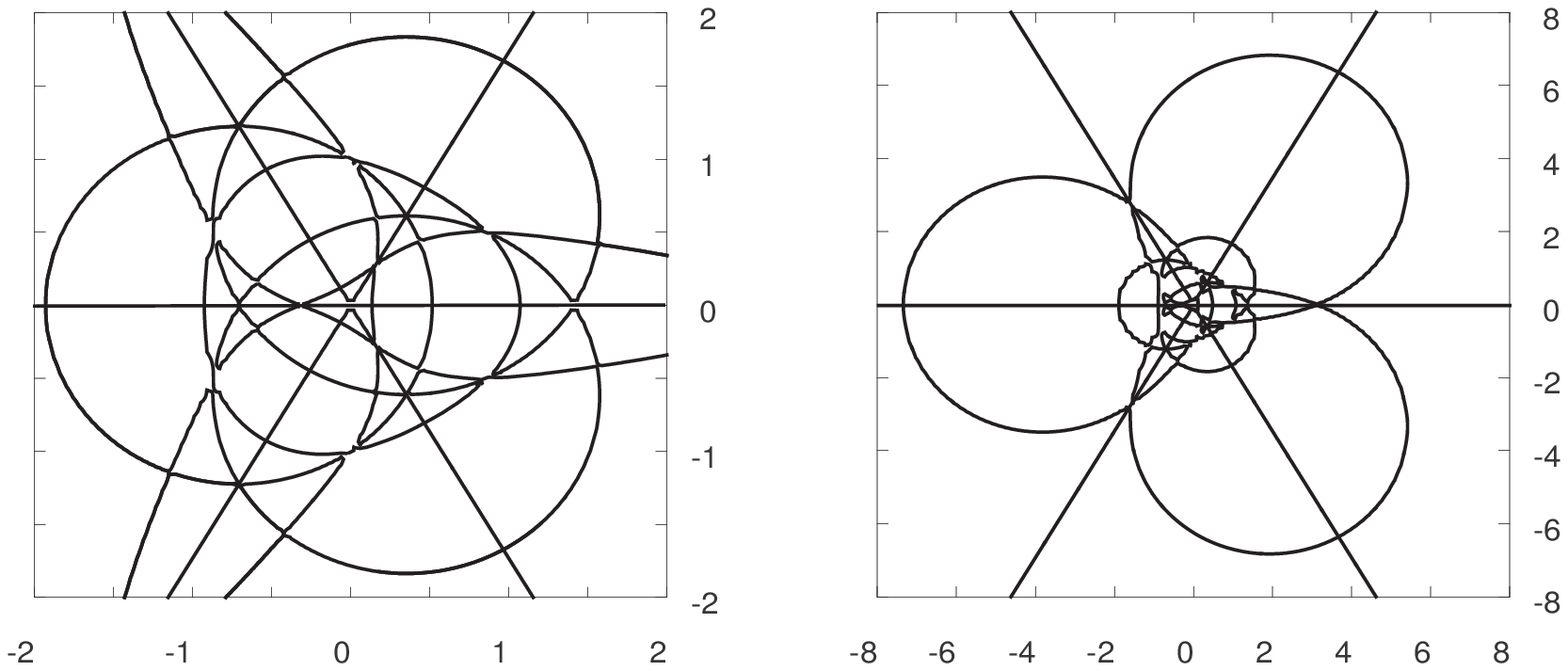}}
\end{center}
\vspace{0.1in}

It is clear from this plot, that the zero set of the symplectic 2-form on the spectral curve projects to a union of 10 
simple closed curves on ${\mathbb{P}}^1$. 

More formally, we can exploit the tetrahedral symmetry to factorize the polynomial $f(\xi,\bar{\xi})$ of 
(\ref{e:lag_curve}).
In our coordinates, the tetrahedral group ${\cal T}$ is the subgroup of SU(2) generated by the 12 fractional linear 
transformations $\xi\mapsto g_k(\xi)e^{v_ji}$ for $j=1,2,3$ and $k=0,1,2,3$, where 
\[
g_0(\xi)=\xi,
\qquad\qquad g_j(\xi)=\frac{\alpha_j\xi-\bar{\beta}_j}{\beta_j\xi+\bar{\alpha}_j},
\]
and 
\[
\alpha=\frac{\sqrt{3}+i}{2\sqrt{3}},
\qquad\qquad \beta_j=\frac{\sqrt{2}i}{\sqrt{3}}e^{-v_ji},
\]
for $v_1=0$, $v_2={\textstyle{\frac{2\pi}{3}}}$ and $v_3={\textstyle{\frac{4\pi}{3}}}$. Note that the spectral curve we
are studying is invariant under this group action, as can be seen by transforming the defining equation (\ref{e:ch3}).

Now, one factor of (\ref{e:lag_curve}) is $\xi-\bar{\xi}$ and so one of the curves is $\xi=\bar{\xi}$. 
Thus one of the Lagrangian curves projects to 
a great circle on ${\mathbb{P}}^1$, and, acting by the tetrahedral group, we  pick out 5 more great circles over which lie 
Lagrangian points. 

In fact, the explicit description of the tetrahedral group above allows us to extract these 6 circles directly, and 
we find the following factor of $f$:
\[
f(\xi,\bar{\xi})=(4\xi^{3} \bar{\xi}^{3}+\sqrt{2} \xi^{3}-18 \xi^{2} \bar{\xi}^{2}+18 \xi \bar{\xi}+\sqrt{2} \bar{\xi}^{3}-4)g(\xi,\bar{\xi}).
\]

Thus the projection of the Lagrangian curves splits naturally into 2 classes: 6 great circles and the solutions
of the remainder of $f$, which is ${\mathbb{R}}{\mbox {e}}(g(\xi,\bar{\xi}))=0$ where
\begin{align}\label{e:factg}
g(\xi,\bar{\xi})=& \xi^{6} \bar{\xi}^{6}+41 \sqrt{2} \xi^{6} \bar{\xi}^{3}-10 \xi^{6}-36 \xi^{5} \bar{\xi}^{5}-9 \sqrt{2} \xi^{5} \bar{\xi}^{2}-126 \xi^{4} \bar{\xi}^{4}+9 \sqrt{2} \xi^{4} \bar{\xi}\nonumber\\
&+302 \xi^{3} \bar{\xi}^{3}-41 \sqrt{2} \xi^{3}-126 \xi^{2} \bar{\xi}^{2}-36 \xi \bar{\xi}+1.
\end{align}

These two classes of curves on ${\mathbb{P}}^1$ are shown below, the left hand being the great circles, while the right 
hand shows the solution set of $\mathbb{R}$e $g(\xi,\bar{\xi})=0$.

\vspace{0.1in}
\setlength{\epsfxsize}{5.0in}
\begin{center}
   \mbox{\epsfbox{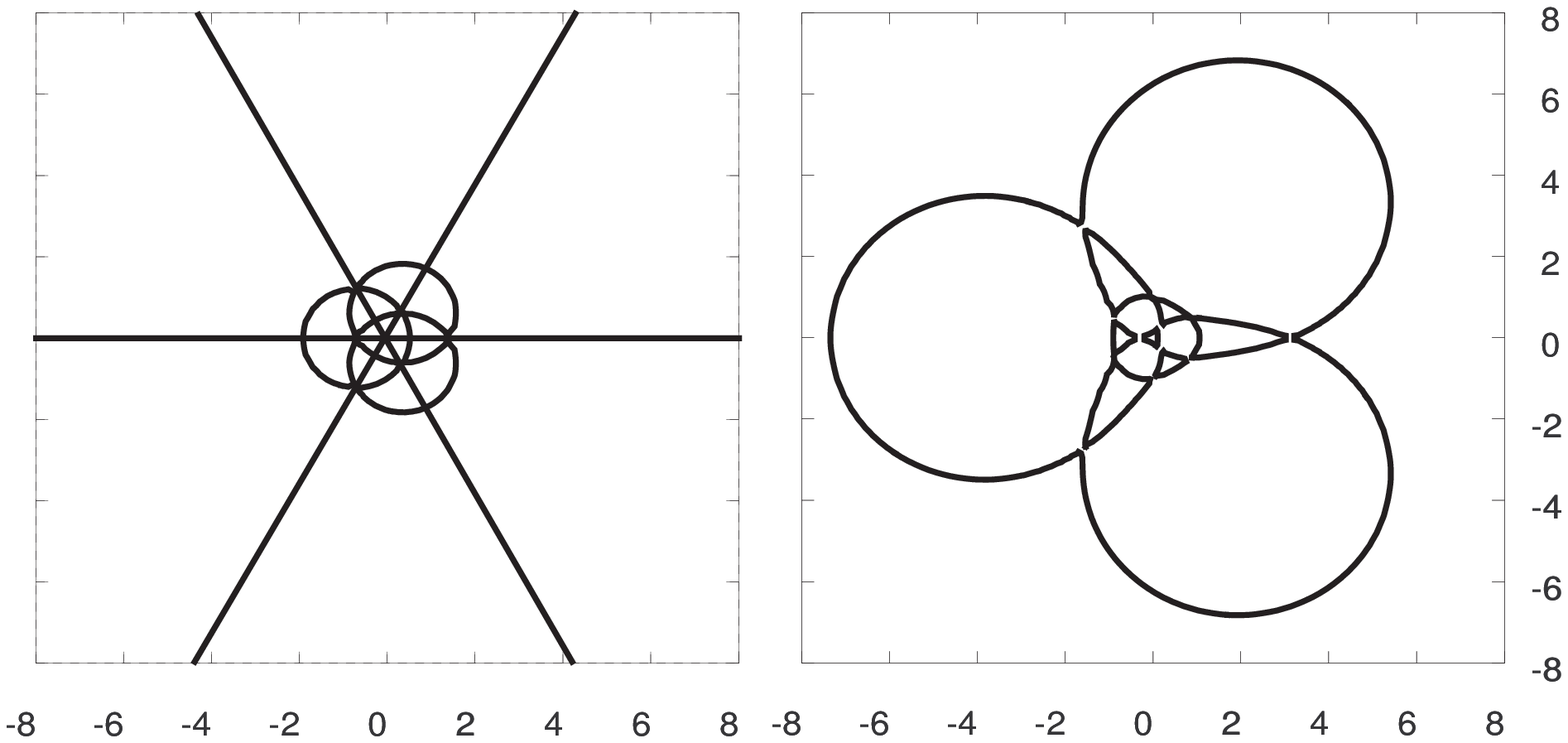}}
\end{center}
\vspace{0.1in}

Let us now look at the ruled surfaces in ${\mathbb{E}}^3$ generated by the great circles. First, these curves contain 
the 6 branch points.  Taking the circle $\xi=\bar{\xi}=s$, and using (\ref{e:coord}) we find the ruled surface to be
\[
x^1=\frac{\Gamma({\textstyle{\frac{1}{3}}})^3(1-s^2)(1-5\sqrt{2}s^3-s^6)^{\scriptstyle{\frac{1}{3}}}
                   +2\sqrt{6}\pi s(1+s^2)r}{\sqrt{6}\pi(1+s^2)^2},
\qquad
x^2=0,
\]
\[
x^3=\frac{-2\Gamma({\textstyle{\frac{1}{3}}})^3s(1-5\sqrt{2}s^3-s^6)^{\scriptstyle{\frac{1}{3}}}
                   +\sqrt{6}\pi (1-s^4)r}{\sqrt{6}\pi(1+s^2)^2}.
\]
The edge of regression of this ruling turns out to be
\[
x^1=\frac{\Gamma({\textstyle{\frac{1}{3}}})^3(1-s^2+s^4)}{\sqrt{6}\pi(1-5\sqrt{2}s^3-s^6)^{\scriptstyle{\frac{2}{3}}}},
\qquad
x^2=0,
\qquad
x^3=-\frac{\Gamma({\textstyle{\frac{1}{3}}})^3s(2-2s^2-5\sqrt{2}s)}{2\sqrt{6}\pi(1-5\sqrt{2}s^3-s^6)^{\scriptstyle{\frac{2}{3}}}}.
\]
Together with the other 5 ruled planes, we get the 6 planes that pass through the edges and centroid of a tetrahedron. 

We turn now to the remaining Lagrangian curves, which project to ${\mathbb{R}}{\mbox {e}}(g(\xi,\bar{\xi}))=0$
with $g$ given by (\ref{e:factg}). While the zero set is made up of 4 simple closed curves, it turns out that
$g$ is not factorizable over ${\mathbb R}$. We argue this as follows.

If $g$ factorized into 4 components, then the tetrahedral action will either leave a component invariant, or move it 
to another component. From the diagram, it is clear that under rotation about the origin through $2\pi/3$ and 
$4\pi/3$ (which is the stereographic projection of a cyclic subgroup of ${\cal T}$) three of the components change
place and one is left invariant. Thus if $g$ factorizes into 4 components, at least one of the components must be invariant
under the cyclic group $C_3$.

Consider now the points where our curve crosses the real axis. If $\xi=u+iv$, then ${\mathbb R}{\mbox e}(g(u,u))$ is 
factorizable over ${\mathbb R}$, in fact, it is equal to:
\[
(u^2-2\sqrt{2}u-1)^2(u^4+5\sqrt{2}u^3-3u^2-5\sqrt{2}u+1)(u^4-\sqrt{2}u^3+3u^2+\sqrt{2}u+1).
\]
Each of these factors (and their products) contains a term that is linear in $u$, and therefore cannot come from
the restriction of a $C_3$-invariant factor to $\xi=u$. Thus, no such factorization exists.

We can still numerically plot the associated (non-planar) ruled surface in ${\mathbb E}^3$, although we cannot 
write down the edge of regression parametrically. 
Combining the information we have obtained on the Lagrangian curves of the charge 3 monopole, we show the 10 
edges of regression in the figure below.

\vspace{0.1in}
\setlength{\epsfxsize}{4.0in}
\begin{center}
   \mbox{\epsfbox{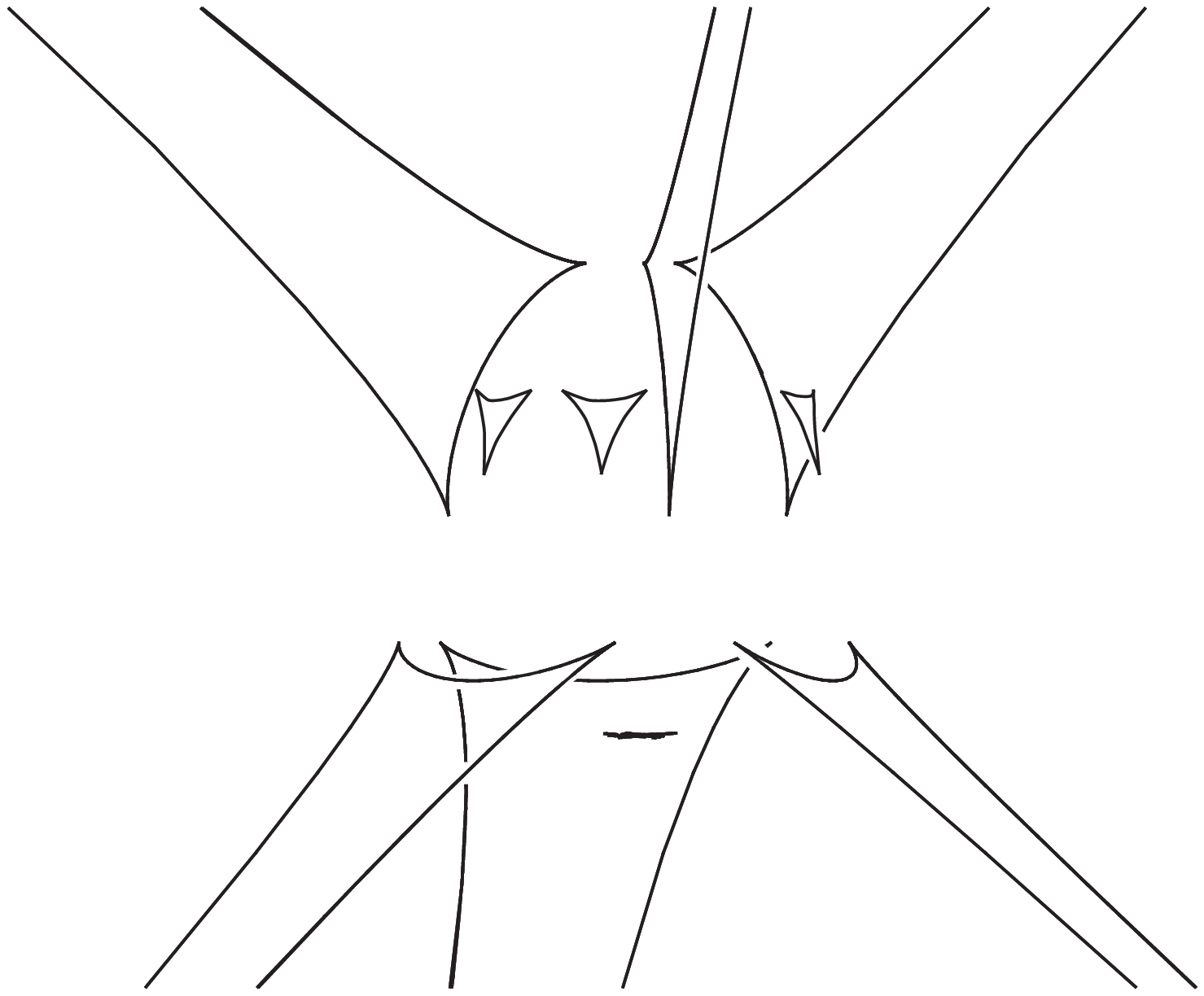}}
\end{center}
\vspace{0.1in}

Note that, once again, the branch points of the spectral curve are the asymptotes of the edges of regression. 

\vspace{0.1in}

\section{Discussion}

The broader geometric context of the preceding is that of K\"ahler surfaces. Thus we consider a real 4-manifold 
${\Bbb{M}}$ endowed
with a complex structure ${\Bbb{J}}$, compatible symplectic structure $\Omega$ and metric ${\Bbb{G}}$. In the
case we have considered, the metric ${\Bbb{G}}$ is of neutral signature (2,2), and so exhibits a rich interplay between 
holomorphic and symplectic structures that is absent for Hermitian metrics. To appreciate this distinction, consider the 
following calibration identity for K\"ahler surfaces. Let $p\in{\Bbb{M}}$ and $v_1,v_2\in \mbox{T}_p{\Bbb{M}}$ span a 
plane. Then  \cite{gak5}
\[
\Omega(v_1,v_2)^2+\epsilon\varsigma^2(v_1,v_2)={\mbox{det }}{\Bbb{G}}(v_i,v_j),
\]
where $\varsigma^2(v_1,v_2)\geq0$ with equality iff $\{v_1,v_2\}$ span a complex plane.
Here, $\epsilon=1$ for ${\Bbb{G}}$ Hermitian, while $\epsilon=-1$ for ${\Bbb{G}}$ neutral. Thus, in the former case 
the above yields a version of the Wirtinger inequality:
\[
{\mbox{det }}{\Bbb{G}}(v_i,v_j)\geq\Omega(v_1,v_2)^2,
\]
while in the latter case we have the metric as a balancing between holomorphic and symplectic structures.
In particular, at a point in a neutral K\"ahler surface a plane can be both holomorphic and Lagrangian - 
a situation that cannot arise in the Hermitian case. 

It is natural then to study complex points on Lagrangian surfaces and Lagrangian points on complex curves in
neutral K\"ahler surfaces. As shown in \cite{gak2}, the complex points on Lagrangian surfaces in T${\mathbb P}^1$ 
correspond precisely to umbilic points on surfaces in  ${\mathbb{E}}^3$. In this paper we considered Lagrangian points 
on holomorphic curves in T${\mathbb P}^1$.

In fact, there is a direct connection between these two situations implicit in the proof of part (iii) of 
Main Theorem 1. We rephrase our result more generally:

\begin{Thm}
Let S be a ruled surface in  ${\mathbb{E}}^3$ associated with a real curve C' in T${\mathbb P}^1$. Then S is flat iff
C' is null.
\end{Thm}
\begin{pf}
The result follows by noting that ${\mathbb {G}}(\dot{\gamma},\dot{\gamma})={\mathbb{I}}{\mbox {m}}\left[(1+\xi\bar{\xi})\dot{\eta}\dot{\bar{\xi}}
                   +2\xi\bar{\eta}\dot{\xi}\dot{\bar{\xi}}\right]$, where
$\dot{\gamma}$ is the tangent vector to the curve C', and recalling equation (\ref{e:gauss}).
\end{pf}

On the one hand, the metric at Lagrangian points on holomorphic curves is zero and so we obtain part (iii) of our Main 
Theorem. On the other hand, the normals along the lines of curvature of a surface S in ${\mathbb L}({\mathbb E}^3)$ are
also null curves in T${\mathbb P}^1$ \cite{gak2}.  Thus we have proven the result in classical surface theory \cite{krey}:

\begin{Thm}
Consider a surface S in ${\mathbb{E}}^3$ and let C be a curve on S. The ruled surface generated by the 1-parameter family 
of normals to S along C is flat iff C is a line of curvature of S.
\end{Thm}

Spectral curves of monopoles have been considered from a number of perspectives. While in principal it is possible to 
reconstruct the
Higgs and Yang-Mills fields in ${\mathbb E}^3$ from the spectral curve via the Nahm data, this procedure is difficult and 
has only been carried out numerically for a small number of symmetric cases (see \cite{ManSut} and references therein). 

One alternative approach has been to consider the minimal surface in ${\mathbb{E}}^3$ generated by the Weierstrass 
representation 
applied to the spectral curve \cite{hitch1}. This has been carried out in detail for the charge 2 case, where the 
resulting geometry has been found to be rich and complicated \cite{small}.
The set of points in ${\mathbb{E}}^3$ where the charge 2 spectral curve lines are orthogonal has also been considered
in \cite{hitch2}, although no explicit calculations were given.
Neither approaches have been extended to higher charge monopoles.

Our techniques, however can be applied to monopoles of any charge. An additional advantage is
that it works for all holomorphic curves and hence avoids the difficulty of the transcendental constraints \cite{HMR}.

Moreover, the symplectic form (and hence the edges of regression) are natural in the following sense.
The Euclidean group O$(3)\ltimes{\mathbb R}^3$ acting on ${\mathbb E}^3$ sends oriented lines to oriented lines and hence
acts on ${\mathbb L}({\mathbb E}^3)$. The symplectic structure $\Omega$, along with ${\mathbb J}$, and hence ${\Bbb{G}}$, 
is invariant under this action.  In fact, up to addition of the round metric on S$^2$, ${\Bbb{G}}$ is the unique metric on 
${\mathbb L}({\mathbb E}^3)$ which is invariant under any subgroup of the Euclidean group \cite{salvai}. 
 
Finally, our techniques can be extended to holomorphic curves in the space of oriented geodesics of any 3-dimensional 
space of constant curvature - in particular to the spectral curves for monopoles in hyperbolic 3-space. 

What is obviously missing is the relationship between the edges of regression and the physical fields (gauge and 
Higgs field). Should this relationship be established, the techniques may yield a new way of localising these fields
in ${\mathbb{E}}^3$ and beyond.

\vspace{0.2in}
\noindent{\bf Acknowledgement}:

The first author would like to thank Nick Manton for suggesting the application to monopoles and Wilhelm Klingenberg
for many helpful discussions, and the second author would like to thank Bernd Kreussler. The second author was supported 
by the IRCSET Embark Initiative Postdoctoral Fellowship Scheme.

\end{document}